\documentclass[12pt]{article}
\usepackage{mathrsfs}
\usepackage{epic,eepic,epsf,epsfig}
\usepackage{amsfonts,srcltx,mathrsfs}
\usepackage{multirow}
\usepackage{amsmath}
\usepackage{amssymb}
\usepackage{amsbsy}
\usepackage{graphicx}
\usepackage{amsfonts}
\usepackage{color}
\usepackage{setspace}
\newtheorem{lem}{Lemma}[section]%
\newtheorem{theorem}[lem]{Theorem}%

\def\nd{\mathrel{\bigm|\kern-.7em/}}

\def\f{\noindent}

\def\P\GammaL{\hbox{\rm P\GammaL}}

\begin{document}
\title{Extremal graphs for disjoint union of vertex-critical graphs}

\footnotetext{E-mails: zhangwq@pku.edu.cn}

\author{Wenqian Zhang\\
{\small School of Mathematics and Statistics, Shandong University of Technology}\\
{\small Zibo, Shandong 255000, P.R. China}}
\date{}
\maketitle

\begin{abstract}
For a graph $F$, let ${\rm EX}(n,F)$ be the set of $F$-free graphs of order $n$ with the maximum number of edges. The graph $F$ is called vertex-critical, if the deletion of its some vertex induces a graph with smaller chromatic number. For example, an odd wheel (obtained by connecting a vertex to a cycle of even length) is a vertex-critical graph with chromatic number 3. For $h\geq2$, let $F_{1},F_{2},...,F_{h}$ be vertex-critical graphs with the same chromatic number. Let $\cup_{1\leq i\leq h}F_{i}$ be the disjoint union of them. In this paper, we characterize the graphs in ${\rm EX}(n,\cup_{1\leq i\leq h}F_{i})$, when there is a proper order among the graphs $F_{1},F_{2},...,F_{h}$. This solves a conjecture (on extremal problem for disjoint union of odd wheels) proposed by Xiao and Zamora \cite{XZ}.

\bigskip

\f {\bf Keywords:} extremal graph; disjoint union of wheels; edge-critical graph; vertex-critical graph.

\end{abstract}

\baselineskip 17 pt

\section{Introduction}

All graphs considered in this paper are finite, undirected and without multi-edges or loops. For a graph $G$, let $\chi(G)$ denote its {\em chromatic number}. The vertex set and edge set of $G$ are denoted by $V(G)$ and $E(G)$, respectively. Let $|G|=|V(G)|$ and $e(G)=|E(G)|$. For a subset $S\subseteq V(G)$, let $G[S]$ be the subgraph of $G$ induced by $S$, and let $G-S=G[V(G)-S]$. For a vertex $u$, denote $G-\left\{u\right\}$ by $G-u$. Another vertex $v$ is called a {\em neighbor} of $u$, if they are adjacent in $G$. Let $N_{G}(u)$ be the set of neighbors of $u$ in $G$, and let $N_{S}(u)=S\cap N_{G}(u)$. Let $d_{G}(u)=|N_{G}(u)|$ and $d_{S}(u)=|N_{S}(u)|$. Let $\delta(G)$ denote the minimum degree of $G$. For any terminology used but not defined here, one may refer to \cite{B}.

 Assume that $G_{1},G_{2},...,G_{\ell}$ are $\ell\geq2$ graphs. Let $\cup_{1\leq i\leq \ell}G_{i}$ be the disjoint union of them, and let $\prod_{1\leq i\leq \ell}G_{i}$  be the graph obtained from $\cup_{1\leq i\leq \ell}G_{i}$ by connecting each vertex in $G_{i}$ to each vertex in $G_{j}$ for any $1\leq i\neq j\leq \ell$. For a positive integer $t$ and a graph $H$, let $tH=\cup_{1\leq i\leq \ell}H$. For a certain integer $n$, let $K_{n}, C_{n}$ and $P_{n}$ denote the complete graph, the cycle and the path of order $n$, respectively. For $r\geq2$, let $K_{n_{1},n_{2},...,n_{r}}$ denote the complete $r$-partite graph with parts of sizes $n_{1},n_{2},...,n_{r}$. Let $T(n,r)$ denote the Tur\'{a}n graph of order $n$ with $r$ parts (i.e., the complete $r$-partite graph of order $n$ with balanced parts).

Let $F$ be a graph. For a graph $G$, denote by $F\subseteq G$ if $G$ contains a copy of $F$, and by $F\nsubseteq G$ otherwise. $G$ is call $F$-free, if $F\nsubseteq G$.   Let ${\rm EX}(n,F)$ be the set of $F$-extremal graphs (i.e., the $F$-free graphs of order $n$ with the maximum number of edges). Let ${\rm ex}(n,F)$ be the number of edges of any graph in ${\rm EX}(n,F)$.
The classical Tur\'{a}n Theorem \cite{T} states that ${\rm EX}(n,K_{r+1})=\left\{T(n,r)\right\}$.  Erd\H{o}s,
Stone and Simonovits \cite{E Simonovits,E Stone} proposed the stability theorem
$${\rm ex}(n,F)=(1-\frac{1}{\chi(F)-1})\frac{n^{2}}{2}+o(n^{2}).$$
The Tur\'{a}n
type problem has attracted many researchers (for example, see \cite{CGPW,DJ,EFGG,HQL,HLF,L,YZ,Y,Y1,Y2}).

A graph $F$ is called {\em edge-critical} if there is an edge $e$ of $F$ such that $\chi(F-e)=\chi(F)-1$, where $F-e$ denotes the graph obtained from $F$ by deleting the edge $e$. A graph $F$ is called {\em vertex-critical} if there is a vertex $u$ of $F$ such that $\chi(F-u)=\chi(F)-1$. Clearly, an edge-critical graph must be vertex-critical. A wheel $W_{n}$ is defined as $W_{n}=K_{1}\prod C_{n-1}$. $W_{n}$ is called an odd wheel if $n$ is odd, and an even wheel otherwise.
Clearly, wheels are vertex-critical graphs, and even wheels are edge-critical graphs. For even wheels, Dzido \cite{D} proved ${\rm ex}(n,W_{2k})=\lfloor\frac{n^{2}}{3}\rfloor$ for $k\geq3$ and $n\geq6k-10$. For odd wheels,
 Dzido and Jastrz\c{e}bski \cite{DJ} determined ${\rm ex}(n,W_{5})$ and ${\rm ex}(n,W_{7})$ for all values of $n$.
Recently, Yuan \cite{Y} determined ${\rm ex}(n,W_{2k+1})$ for $k\geq3$ and sufficiently large $n$.

A nearly $k$-regular graph is a graph of which all vertices have degree $k$ except one vertex with degree $k-1$. Let $\mathcal{U}^{k-1}_{n}(P_{2k-1})$ be
the set of $P_{2k-1}$-free, $(k-1)$-regular or nearly $(k-1)$-regular graphs of order $n$. Let $\mathcal{K}_{n_{1},n_{2}}(\mathcal{U}^{k-1}_{n_{1}}(P_{2k-1});P_{2})$ be the set of graphs obtained from $K_{n_{1},n_{2}}$ by embedding a graph from $\mathcal{U}^{k-1}_{n_{1}}(P_{2k-1})$ in the part of size $n_{1}$, and embedding an edge in the part of size $n_{2}$.

\begin{theorem}{\rm (Yuan \cite{Y})}\label{1}
For $k\geq3$ and sufficiently large $n$,
$${\rm ex}(n, W_{2k+1})=\max_{1\leq n_{0}\leq n}\left\{n_{0}(n-n_{0})+\left\lfloor\frac{k-1}{2}n_{0}\right\rfloor+1\right\}.
$$
 Moreover, ${\rm EX}(n,W_{2k+1})\subseteq\mathcal{K}_{n_{0},n-n_{0}}(\mathcal{U}^{k-1}_{n_{0}}(P_{2k-1});P_{2})$.
\end{theorem}

The extremal problem for disjoint union of graphs has been paid attention (for example, see (\cite{G,LLP,XZ,YZ})). In particular, Xiao and Zamora \cite{XZ} studied the Tur\'{a}n number of disjoint union of wheels. Since even wheels are edge-critical, for large $n$, the Tur\'{a}n number of disjoint union of even wheels can be easily determined by a result of Simonovits \cite{S}. For the disjoint union of odd wheels, they obtained the following result.

\begin{theorem}{\rm (Xiao and Zamora \cite{XZ})}\label{2}
For $k\geq3,m\geq2$ and sufficiently large $n$,
\begin{equation}
\begin{aligned}
&{\rm ex}(n,mW_{2k+1})\\
&=\max_{1\leq n_{0}\leq n}\left\{\tbinom{m-1}{2}+\left\lfloor\frac{k-1}{2}n_{0}\right\rfloor+(n_{0}+m-1)(n-m+1)-n^{2}_{0}+1\right\}\\
&=\tbinom{m-1}{2}+(m-1)(n-m+1)+{\rm ex}(n-m+1,W_{2k+1}).
\end{aligned}\notag
\end{equation}
 Moreover, each graph in ${\rm EX}(n,mW_{2k+1})$ is of form $K_{m-1}\prod H$, where $$H\in\mathcal{K}_{n_{0},n-m+1-n_{0}}(\mathcal{U}^{k-1}_{n_{0}}(P_{2k-1});P_{2}).$$
\end{theorem}

For the disjoint union of distinct odd wheels,  Xiao and Zamora \cite{XZ} proposed the following conjecture (which follows from Theorem \ref{main}).

\medskip

\f{\bf Conjecture 1} {\rm (Xiao and Zamora \cite{XZ})} Assume that $k_{1}\geq k_{2}\geq\cdots\geq k_{m}$. For sufficiently large $n$,
 \begin{equation}
\begin{aligned}
&{\rm ex}(n,\cup_{1\leq i\leq m} W_{2k_{i}+1})\\
&=\max_{1\leq i\leq m,1\leq n_{0}\leq n}\left\{n_{0}(n-n_{0})+(i-1)(n_{0}-i+1)+\tbinom{i-1}{2}+\left\lfloor\frac{k_{i}-1}{2}(n_{0}-i+1)\right\rfloor+1\right\}\\
&=\max_{1\leq i\leq m}\left\{\tbinom{i-1}{2}+(i-1)(n-i+1)+{\rm ex}(n-i+1,W_{2k_{i}+1})\right\}.
\end{aligned}\notag
\end{equation}

\medskip

Let $F_{1},F_{2},...,F_{h}$ be $h\geq1$ graphs with $\chi(F_{i})=r+1\geq3$ for any $1\leq i\leq h$. We say that $F_{1},F_{2},...,F_{h}$ is properly ordered, if for any $1\leq \ell\leq h$, there exists at least one graph in ${\rm EX}(n,F_{\ell})$ which contains no copy of any $F_{j}$ with $1\leq j\leq\ell$. The main result of this paper is the following theorem.

\begin{theorem}\label{main}
For $h\geq2$, let $F_{1},F_{2},...,F_{h}$ be properly ordered vertex-critical graphs with $\chi(F_{i})=r+1\geq3$ for any $1\leq i\leq h$. Then, for sufficiently large $n$,
$${\rm ex}(n,\cup_{1\leq i\leq h}F_{i})=\max_{1\leq\ell\leq h}\left\{\tbinom{\ell-1}{2}+(\ell-1)(n-\ell+1)+{\rm ex}(n-\ell+1,F_{\ell})\right\}.$$
 Moreover, each graph in ${\rm EX}(n,\cup_{1\leq i\leq h}F_{i})$ is of form $K_{\ell-1}\prod H$, where $H\in{\rm EX}(n-\ell+1,F_{\ell})$.
\end{theorem}

For any $k\geq2$ and sufficiently large $n$, from Theorem \ref{1} we see that each graph in ${\rm EX}(n,W_{2k+1})$ contains no copy of any $W_{2\ell+1}$ with $\ell\geq k$. Thus, $W_{2k_{1}+1},W_{2k_{2}+1},...,W_{2k_{h}+1}$ is properly ordered if $k_{1}\geq k_{2}\geq\cdots\geq k_{h}$.
Recall that odd wheels are vertex-critical graphs. Hence, Conjecture 1 follows directly  from Theorem \ref{main}.

Let $F$ be a vertex-critical graph with $\chi(F)\geq3$. Clearly, $F,F,...,F$ is properly ordered by definition. Thus, the following theorem follows directly from Theorem \ref{main}.

\begin{theorem}
Let $F$ be a vertex-critical graph with $\chi(F)\geq3$. For $h\geq2$ and sufficiently large $n$,
$${\rm ex}(n,hF)=\max_{1\leq\ell\leq h}\left\{\tbinom{\ell-1}{2}+(\ell-1)(n-\ell+1)+{\rm ex}(n-\ell+1,F)\right\}$$
$$=\tbinom{h-1}{2}+(h-1)(n-h+1)+{\rm ex}(n-h+1,F).$$
 Moreover, each graph in ${\rm EX}(n,hF)$ is of form $K_{h-1}\prod H$, where $H\in{\rm EX}(n-h+1,F)$.
\end{theorem}

The rest of this paper is organized as follows. In Section 2, we present several well-known results which will be used in the proof of the main result. In Section 3, we give the proof of Theorem \ref{main}.

\section{Preliminaries}

 The following result is the classical stability theorem proved by  Erd\H{o}s,
Stone and Simonovits \cite{E Simonovits,E Stone}:

\begin{theorem} {\rm (Erd\H{o}s, Stone and Simonovits \cite{E Simonovits,E Stone})}\label{stability}
Let $F$ be a graph with $\chi(F)=r+1\geq3$. For every $\epsilon>0$, there exist $\delta>0$ and $n_{0}$ such that, if $G$ is an $F$-free graph of order $n\geq n_{0}$ with $e(G)>(\frac{r-1}{2r}-\delta)n^{2}$, then $G$ can be
obtained from $T(n,r)$ by adding and deleting at most $\epsilon n^{2}$ edges.
\end{theorem}

The following well-known lemma is proved by Erd\H{o}s \cite{E} and
Simonovits \cite{S}.

\begin{lem}{\rm (Erd\H{o}s \cite{E} and Simonovits \cite{S})}\label{mini degree}
Let $F$ be a graph with $\chi(F)=r+1\geq3$. For any given $\theta>0$, if $n$ is sufficiently large, then $\delta(G)>(1-\frac{1}{r}-\theta)n$ for any $G\in{\rm EX}(n,F)$.
\end{lem}

\section{Proof of Theorem \ref{main}}

The following fact will be used.

\medskip

\f{\bf Fact 1.}
Let $A_{1},A_{2},...,A_{\ell}$ be $\ell\geq2$ finite sets. Then
$$|\cap_{1\leq i\leq \ell}A_{i}|\geq(\sum_{1\leq i\leq\ell}|A_{i}|)-(\ell-1)|\cup_{1\leq i\leq \ell}A_{i}|.$$

\medskip

Now we can give the proof of Theorem \ref{main}.

\medskip

\f{\bf Proof of Theorem \ref{main}.} Recall that $F_{i}$ is a vertex-critical graph with $\chi(F_{i})=r+1\geq3$ for any $1\leq i\leq h$. Moreover, $F_{1},F_{2},...,F_{h}$ is properly ordered. Thus, for any $1\leq \ell\leq h$, we can choose a graph $H_{\ell}$  in ${\rm EX}(n-\ell+1,F_{\ell})$ such that $H_{\ell}$ contains no copy of any $F_{i}$ with $1\leq i\leq \ell$. Let $G(\ell)=K_{\ell-1}\prod H_{\ell}$. Now we show that $G(\ell)$ contains no copy of $\cup_{1\leq i\leq h}F_{i}$. Suppose  not. Let $x_{1},x_{2},...,x_{\ell-1}$ be the vertices of $G(\ell)$ contained in the part $K_{\ell-1}$. Then $G(\ell)-\left\{x_{1},x_{2},...,x_{\ell-1}\right\}=H_{\ell}$ contains a copy of $F_{j}$ for some $1\leq j\leq \ell$, a contradiction. Hence, $G(\ell)$ contains no copy of $\cup_{1\leq i\leq h}F_{i}$ for any $1\leq \ell\leq r$.

Set $t=\sum_{1\leq i\leq h}|F_{i}|$. Let $G\in {\rm EX}(n,\cup_{1\leq i\leq h}F_{i})$, where $n$ is sufficiently large. Then $e(G)\geq \max_{1\leq \ell\leq h}e(G(\ell))$.
Since $T(n,r)$ is $\cup_{1\leq i\leq h}F_{i}$-free, we have $e(G)\geq e(T(n,r))$.  As is well known, $e(T(n,r))\geq\frac{r-1}{2r}n^{2}-\frac{r}{8}$. Thus, the following Claim 1 holds.

\medskip

\f{\bf Claim 1.} $e(G)\geq\frac{r-1}{2r}n^{2}-\frac{r}{8}$.

\medskip

Let $\theta>0$ be a fixed constant with respect to $n$, which is small enough for all the inequalities appeared in the following.

\medskip

\f{\bf Claim 2.}
There is a partition $V(G)=\cup_{1\leq i\leq r}V_{i}$ such that $\sum_{1\leq i\leq r}e(G[V_{i}])$ is minimum, where $\sum_{1\leq i\leq r}e(G[V_{i}])\leq\theta^{3}n^{2}$ and $||V_{i}|-\frac{n}{r}|\leq\theta n$ for any $1\leq i\leq r$.

\medskip

\f{\bf Proof of Claim 2.} Since $e(G)\geq\frac{r-1}{2r}n^{2}-\frac{r}{8}$ from Claim 1, by Theorem \ref{stability} (letting $\epsilon=\theta^{3}$ and $F=\cup_{1\leq i\leq h}F_{i}$ there), $G$ is obtained from $T(n,r)$ by deleting and adding at most $\theta^{3}n^{2}$ edges. That is to say, there is a balanced partition $V(G)=\cup_{1\leq i\leq r}U_{i}$ such that $\sum_{1\leq i\leq r}e(G[U_{i}])\leq\theta^{3}n^{2}$. Let $V(G)=\cup_{1\leq i\leq r}V_{i}$ be a partition such that $\sum_{1\leq i\leq r}e(G[V_{i}])$ is minimum. Then
 $$\sum_{1\leq i\leq r}e(G[V_{i}])\leq\sum_{1\leq i\leq r}e(G[U_{i}])\leq\theta^{3}n^{2}.$$
Let $a=\max_{1\leq i\leq r}||V_{i}|-\frac{n}{r}|$. Without loss of generality, assume that $a=||V_{1}|-\frac{n}{r}|$. Using the Cauchy-Schwarz inequality, we obtain that
$$2\sum_{2\leq i<j\leq r}|V_{i}||V_{j}|=(\sum_{2\leq i\leq r}|V_{i}|)^{2}-\sum_{2\leq i\leq r}|V_{i}|^{2}\leq\frac{r-2}{r-1}(\sum_{2\leq i\leq r}|V_{i}|)^{2}=\frac{r-2}{r-1}(n-|V_{1}|)^{2}.$$
Thus,
\begin{equation}
\begin{aligned}
e(G)&\leq(\sum_{1\leq i<j\leq r}|V_{i}||V_{j}|)+(\sum_{1\leq i\leq r}e(G[V_{i}]))\\
&\leq|V_{1}|(n-|V_{1}|)+(\sum_{2\leq i<j\leq r}|V_{i}||V_{j}|)+\theta^{3}n^{2}\\
&\leq|V_{1}|(n-|V_{1}|)+\frac{r-2}{2(r-1)}(n-|V_{1}|)^{2}+\theta^{3}n^{2}\\
&=\frac{r-1}{2r}n^{2}-\frac{r}{2(r-1)}a^{2}+\theta^{3}n^{2}.
\end{aligned}\notag
\end{equation}
Recall that $e(G)\geq\frac{r-1}{2r}n^{2}-\frac{r}{8}\geq(\frac{r-1}{2r}-\theta^{3})n^{2}$. It follows that $a\leq\sqrt{\frac{4(r-1)}{r}\theta^{3}n^{2}}\leq\theta n$. This finishes the proof of Claim 2. \hfill$\Box$

\medskip

For $1\leq i\leq r$, let $W_{i}=\left\{v\in V_{i}~|~d_{V_{i}}(v)\geq\theta n\right\}$,  and let $W=\cup_{1\leq i\leq r}W_{i}$.

\medskip

\f{\bf Claim 3.}
$|W|\leq 2\theta^{2}n$.

\medskip

\f{\bf Proof of Claim 3.}  Since $\sum_{1\leq i\leq r}e(G[V_{i}])\leq\theta^{3}n^{2}$, and
$$\sum_{1\leq i\leq r}e(G[V_{i}])=\sum_{1\leq i\leq r}\frac{1}{2}\sum_{v\in V_{i}}d_{V_{i}}(v)\geq\sum_{1\leq i\leq r}\frac{1}{2}\sum_{v\in W_{i}}d_{V_{i}}(v)\geq\frac{1}{2}\sum_{1\leq i\leq r}\theta n|W_{i}|=\frac{1}{2}\theta n|W|,$$
we have $|W|\leq2\theta^{2}n$.
This finishes the proof of Claim 3.
\hfill$\Box$

\medskip

For any $1\leq i\leq r$, let $\overline{V}_{i}=V_{i}-W$.

\medskip

\f{\bf Claim 4.}
Let $1\leq \ell\leq r$ be fixed. Assume that  $u_{1},u_{2},...,u_{rt}\in \cup_{1\leq i\neq \ell\leq r}\overline{V}_{i}$ and $u\in W_{i_{0}}$ with $i_{0}\neq \ell$. Then there are $3t$ vertices in $\overline{V}_{\ell}$ which are adjacent to all the vertices $u,u_{1},u_{2},...,u_{rt}$ in $G$.

\medskip

\f{\bf Proof of Claim 4.} Recall that $\frac{n}{r}-\theta n\leq|V_{s}|\leq\frac{n}{r}+\theta n$ for any $1\leq s\leq r$ by Claim 2. By Claim 3, we have $\frac{n}{r}-2\theta n\leq|\overline{V}_{s}|\leq\frac{n}{r}+\theta n$. Since $n$ is sufficiently large, we have $\delta(G)\geq(\frac{r-1}{r}-\theta)n$ by Lemma \ref{mini degree}. Since $V(G)=\cup_{1\leq i\leq r}V_{i}$ is a partition such that $\sum_{1\leq i\leq r}e(G[V_{i}])$ is minimum, we have $d_{V_{\ell}}(u)\geq d_{V_{i_{0}}}(u)$ as $u\in V_{i_{0}}$.  Otherwise, we will obtain a contradiction by moving $u$ from $V_{i_{0}}$ to $V_{\ell}$. It follows that
$$d_{V_{\ell}}(u)\geq\frac{d_{V_{\ell}}(u)+d_{V_{i_{0}}}(u)}{2}\geq\frac{d_{G}(u)-\sum_{1\leq j\neq\ell,i_{0}\leq r}|V_{j}|}{2}\geq(\frac{1}{2r}-(r-1)\theta)n.$$
Then, using Claim 3,
 $$d_{\overline{V}_{\ell}}(u)\geq d_{V_{\ell}}(u)-|W|\geq(\frac{1}{2r}-r\theta)n.$$

 For any $1\leq i\leq rt$, assume that $u_{i}\in \overline{V}_{j_{i}}$, where $j_{i}\neq\ell$. Then $d_{V_{j_{i}}}(u_{i})\leq\theta n$ as $u_{i}\notin W$.
Hence $$d_{V_{\ell}}(u_{i})\geq d_{G}(u_{i})-d_{V_{j_{i}}}(u_{i})-\sum_{1\leq s\neq\ell,j_{i}\leq r}|V_{s}|\geq(\frac{1}{r}-r\theta)n.$$
Thus,
 $$d_{\overline{V}_{\ell}}(u_{i})\geq d_{V_{\ell}}(u_{i})-|W|\geq(\frac{1}{r}-2r\theta)n.$$
By Fact 1, we have
$$|N_{\overline{V}_{\ell}}(u)\cap(\cap_{1\leq i\leq rt}N_{\overline{V}_{\ell}}(u_{i}))|\geq|N_{\overline{V}_{\ell}}(u)|+(\sum_{1\leq i\leq rt}|N_{\overline{V}_{\ell}}(u_{i})|)-rt|V_{\ell}|\geq(\frac{1}{2r}-rt(r+2)\theta)n\geq 3t.$$
Thus, there are $3t$ vertices in $\overline{V}_{\ell}$ which are adjacent to all the vertices $u,u_{1},u_{2},...,u_{rt}$.
 This finishes the proof of Claim 4. \hfill$\Box$

\medskip

\f{\bf Claim 5.}
$|W|\leq h-1$. Moreover, for any $S\subseteq V(G)-W$ with $|S|\leq t$, $G$ contains the disjoint union of $F'_{1},F'_{2},...,F'_{|W|}$, where $F'_{i}$ is a copy of $F_{i}$, and $|W\cap V(F'_{i})|=1$ and $V(F'_{i})\cap S=\emptyset$ for any $1\leq i\leq |W|$.

\medskip

\f{\bf Proof of Claim 5.}  Suppose that $|W|\geq h$. We will prove that $G$ contains the disjoint union of $F'_{1},F'_{2},...,F'_{h}$, where $F'_{i}$ is a copy of $F_{i}$, and $|W\cap V(F'_{i})|=1$ and $V(F'_{i})\cap S=\emptyset$ for any $1\leq i\leq h$. Assume that we have find such $F'_{1},F'_{2},...,F'_{h_{0}}$, where $h_{0}\leq h-1$. Recall that $|W|\geq h$ and $|W\cap V(F'_{i})|=1$ for any $1\leq i\leq h_{0}$. We can find $w\in W-S\cup(\cup_{1\leq i\leq h_{0}}V(F'_{i}))$.  Without loss of generality, assume that $w\in W_{1}$. Then, using Claim 3, $d_{\overline{V}_{1}}(w)\geq d_{V_{1}}(w)-|W|\geq\frac{1}{2}\theta n\geq 3t$. Thus, we can find $t$ neighbors of $w$ in $\overline{V}_{1}-S\cup(\cup_{1\leq i\leq h_{0}}V(F'_{i}))$, say $w^{1}_{1},w^{1}_{2},...,w^{1}_{t}$. By Claim 4, there are $3t$ vertices in $\overline{V}_{2}$ which are adjacent to all the vertices $w,w^{1}_{1},w^{1}_{2},...,w^{1}_{t}$. Among these $3t$ vertices, there are $t$ ones in  $\overline{V}_{2}-S\cup(\cup_{1\leq i\leq h_{0}}V(F'_{i}))$, say $w^{2}_{1},w^{2}_{2},...,w^{2}_{t}$. Repeat this process using Claim 4. We can obtain vertices $w^{j}_{1},w^{j}_{2},...,w^{j}_{t}$ in  $\overline{V}_{j}-S\cup(\cup_{1\leq i\leq h_{0}}V(F'_{i}))$ for any $1\leq j\leq r$, which are adjacent to all the vertices $w,w^{1}_{1},w^{1}_{2},...,w^{j-1}_{t}$. Thus $G-S\cup(\cup_{1\leq i\leq h_{0}}V(F'_{i}))$ contains a copy of $K_{1}\prod T(rt,r)$. Note that $F_{h_{0}+1}\subseteq K_{1}\prod T(rt,r)$, since $F_{h_{0}+1}$ is vertex-critical. Hence, $G$ contains the disjoint union of $F'_{1},F'_{2},...,F'_{h}$, where $F'_{i}$ is a copy of $F_{i}$, and $|W\cap V(F'_{i})|=1$ and $V(F'_{i})\cap S=\emptyset$ for any $1\leq i\leq h$. This is impossible, since $G$ is $\cup_{1\leq i\leq h}F_{i}$-free. Thus, $|W|\leq h-1$.
The desired disjoint union of $F'_{1},F'_{2},...,F'_{|W|}$ can be found similarly.
This finishes the proof of Claim 5. \hfill$\Box$

\medskip

Set $|W|=q$. Then $q\leq h-1$.

\medskip

\f{\bf Claim 6.}
$G-W$ contains no copy of $F_{q+1}$.

\medskip

\f{\bf Proof of Claim 6.}
Suppose not. $G-W$ contains a $F'_{q+1}$, where $F'_{q+1}$ is a copy of $F_{q+1}$. If $G-W$ contains a copy of $\cup_{q+1\leq j\leq h}F_{j}$, then $G$ contain a copy of  $\cup_{1\leq j\leq h}F_{j}$ by Claim 5 (letting $S=\cup_{q+1\leq j\leq h}V(F_{j})$), a contradiction. Thus, we can assume that $G-W$ contains a $\cup_{q+1\leq j\leq m}F'_{j}$ for some $q+1\leq m\leq h-1$, where $F'_{j}$ is a copy of $F_{j}$ for any $q+1\leq j\leq m$. Moreover, $G-W\cup (\cup_{q+1\leq j\leq m}V(F'_{j}))$ contains no copy of $F_{m+1}$. Let $S=\cup_{q+1\leq j\leq m}V(F'_{j})$. Then $|S|\leq t$, and $G-W\cup S$ contains no copy of $F_{m+1}$.

Now we show that, (if exists) any subgraph $F'_{m+1}$ of $G-W$ as a copy of $F_{m+1}$, must contain an edge inside $S$ or an edge connecting $S\cap \overline{V}_{i}$ to $\overline{V}_{i}-S$ for some $1\leq i\leq r$. Suppose not for $F'_{m+1}$. Let $S_{i}=V(F'_{m+1})\cap S\cap \overline{V}_{i}$ for any $1\leq i\leq r$. Then $S\cap V(F'_{m+1})=\cup_{1\leq i\leq r}S_{i}$ is an independent set of $F'_{m+1}$, and for any $1\leq \ell\leq r$, the vertices in $S_{\ell}$ are only adjacent to the vertices in $V(F'_{m+1})\cap(\cup_{1\leq i\neq\ell\leq r}(\overline{V}_{i}-S_{i}))$ in the graph $F'_{m+1}$. Note that $|S|+|F'_{m+1}|\leq 2t$. By Claim 4, there are $t$ vertices in $\overline{V}_{\ell}-S\cup V(F'_{m+1})$, say $v^{\ell}_{1},...,v^{\ell}_{t}$, which  are adjacent to all the vertices in $V(F'_{m+1})\cap(\cup_{1\leq i\neq\ell\leq r}(\overline{V}_{i}-S_{i}))$. Let $F''_{m+1}$ be the induced subgraph of $G-W\cup S$ obtained from $F'_{m+1}$ by deleting the vertices in $S\cap V(F'_{m+1})$, and adding the vertices $v^{\ell}_{1},...,v^{\ell}_{t}$ for all $1\leq \ell\leq r$. Then $F''_{m+1}$ contains a copy of $F'_{m+1}$ (or $F_{m+1}$). That is to say, $G-W\cup S$ contains a copy of $F_{m+1}$, a contradiction. Therefore, any subgraph $F'_{m+1}$ of $G-W$ as a copy of $F_{m+1}$, must contain an edge inside $S$ or an edge connecting $S\cap \overline{V}_{i}$ to $\overline{V}_{i}-S$ for some $1\leq i\leq r$.

Choose a vertex $u$ in $S=\cup_{q+1\leq j\leq m}V(F'_{j})$. Without loss of generality, assume $u\in \overline{V}_{1}$. Let $G'$ be the graph obtained from $G$ by deleting all the edges inside $S$ and edges connecting $S\cap \overline{V}_{i}$ to $\overline{V}_{i}-S$ for all $1\leq i\leq r$, and adding all the edges between $u$ and $\overline{V}_{1}-S$. Clearly, we delete at most $\tbinom{|S|}{2}+|S|\theta n\leq t^{2}+t\theta n\leq2t\theta n$ edges, and add at least $(\frac{1}{r}-\theta)n-|W|-|S|\geq(\frac{1}{r}-2\theta)n$ edges. It follows that
$$e(G')-e(G)\geq (\frac{1}{r}-2\theta)n-2t\theta n\geq(\frac{1}{r}-3t\theta)n>0.$$
 Thus, $e(G')>e(G)\geq e(G(m+1))$. By the discussion of the last paragraph, we see that $G'-W\cup \left\{u\right\}$ contains no copy of $F_{m+1}$. Noting $m\geq q+1$, any induced subgraph of $G'-W\cup \left\{u\right\}$ with $n-m$ vertices has at most ${\rm ex}(n-m,F_{m+1})$ edges.  Recall that $G(m+1)=K_{m}\prod H_{m}$ and $e(H_{m})={\rm ex}(n-m,F_{m+1})$. It follows that $e(G')\leq e(G(m+1))$, a contradiction. Therefore, $G-W$ contains no copy of $F_{q+1}$. This finishes the proof of Claim 6. \hfill$\Box$

\medskip

Since $G-W$ contains no copy of $F_{|W|+1}$ by Claim 6, we have $e(G-W)\leq {\rm ex}(n-|W|,F_{|W|+1})$. It follows that $e(G)\leq G(|W|+1)$. Recall that $e(G)\geq \max_{1\leq \ell\leq r}e(G(\ell))$. Consequently, $e(G)= G(|W|+1)=\max_{1\leq \ell\leq r}e(G(\ell))$, and $G=K_{|W|}\prod H$, where $H\in{\rm EX}(n-|W|,F_{|W|+1})$. This completes the proof. \hfill$\Box$

\medskip

\f{\bf Declaration of competing interest}

\medskip

There is no conflict of interest.

\medskip

\f{\bf Data availability statement}

\medskip

No data was used for the research described in the article.


\medskip

\begin{thebibliography}{99}






\bibitem{B}
 B. Bollob\'{a}s, Extremal Graph Theory, Academic Press, New York, 1978.





\bibitem{CGPW}
 G. Chen, R. Gould, F. Pfender and B. Wei, Extremal graphs for intersecting cliques, J. Comb. Theory,
Ser. B 89 (2003) 159-171.

\bibitem{D}
 T. Dzido, A note on Tur\'{a}n numbers for even wheels, Graphs Comb. 29 (2013) 1305-1309.

\bibitem{DJ}
 T. Dzido and A. Jastrz\c{e}bski, Tur\'{a}n numbers for odd wheels, Discrete Math. 341 (4) (2018) 1150-1154.



\bibitem{E}
 P. Erd\H{o}s, On some new inequalities concerning extremal properties of graphs, stability problems, Theory Graphs
(Proc. Colloq. Tihany, 1966) (1968) 77-81.

\bibitem{EFGG}
 P. Erd\H{o}s, Z. F\"{u}redi, R. Gould and D. Gunderson, Extremal graphs for intersecting triangles, J.
Comb. Theory, Ser. B 64 (1995) 89-100.


\bibitem{E Simonovits}
 P. Erd\H{o}s and M. Simonovits, A limit theorem in graph theory, Studia Sci. Math. Hungar. 1 (1966) 51-57.

\bibitem{E Stone}
 P. Erd\H{o}s and A. Stone, On the structure of linear graphs, Bull. Am. Math. Soc. 52 (1946) 1087-1091.

\bibitem{G}
 I. Gorgol, Tur\'{a}n numbers for disjoint copies of graphs, Graphs Comb. 27 (2011) 661-667.

\bibitem{HQL}
 X. Hou, Y. Qiu and B. Liu, Extremal graph for intersecting odd cycles, Electron. J. Comb. 23 (2) (2016) P2.29.


\bibitem{HLF}
X. He, Y. Li and L. Feng, Extremal graphs for the odd prism, Discrete Math. 348 (2025) 114249.

\bibitem{L}
H. Liu, Extremal graphs for blow-ups of cycles and trees, Electron. J. Comb. 20 (1) (2013) P65.

\bibitem{LLP}
 B. Lidick\'{y}, H. Liu and C. Palmer, On the Tur\'{a}n number of forests, Electron. J. Comb. 20 (2) (2013) P62.


\bibitem{S}
 M. Simonovits, A method for solving extremal problems in graph theory, stability problems, in: Theory of Graphs
(Proc. Colloq., Tihany, 1966), 1968, pp. 279-319.



\bibitem{T}
 P. Tur\'{a}n, On an extremal problem in graph theory, Mat. Fiz. Lapok 48 (1941) 436-452.

\bibitem{XZ}
C. Xiao and O. Zamora, A note on the Tur\'{a}n number of disjoint union of wheels, Discrete Math. 344 (2021) 112570.

\bibitem{YZ}
 L. Yuan, and X. Zhang, Tur\'{a}n numbers for disjoint paths, J. Graph Theory 98 (2021) 499-524.
\bibitem{Y}
L. Yuan, Extremal graphs for wheels, J. Graph Theory 98 (2021) 691-707.
\bibitem{Y1}
 L. Yuan, Extremal graphs for edge-blow up of graphs, J. Comb. Theory, Ser. B 152 (2022) 379-398.
\bibitem{Y2}
 L. Yuan, Extremal graphs for the $k$-flower, J. Graph Theory 89 (2018) 26-39.






\end{thebibliography}
\end{document}